\documentclass[11pt]{article}
\usepackage{amssymb}
\usepackage{latexsym,bm}
\usepackage{graphicx}
\usepackage{mathrsfs,amscd,amssymb,amsthm,amsmath,bm,graphicx,psfrag,subfigure,url}
\usepackage{amsmath}
\usepackage{mathrsfs}

\date{}
\begin{document}
\title{ The minimum Kirchhoff index of phenylene chains \footnote{E-mail addresses:
{\tt mathdzhang@163.com}(L.Zhang).}}
\author{\hskip -10mm  Leilei Zhang\thanks{Corresponding author.}\\
{\hskip -10mm \small  Department of Mathematics, East China Normal University, Shanghai 200241, China}}\maketitle
\begin{abstract}
Let $G$ be a connected graph. The resistance distance between any two vertices of $G$ is equal to the effective resistance between them in the corresponding electrical network constructed from $G$ by replacing each edge with a unit resistor. The Kirchhoff index is defined as the sum of resistance distances between all pairs of the vertices. Recently, Yang and Wang determined the maximum Kirchhoff index of phenylene chains, and they proposed a conjecture about the minimum Kirchhoff index. In this note, we characterized the minimum phenylene chains with respect to the Kirchhoff index. This proves the conjecture.
\end{abstract}

{\bf Key words.} Resistance distance; Kirchhoff index; phenylene chains; $S, T$-isomers

{\bf Mathematics Subject Classification.} 05C09, 05C92, 05C12
\vskip 8mm

\section{ Introduction}
In this paper, we consider simple connected graph $G=(V(G),E(G))$ with vertex set $V(G)$ and edge set $E(G)$. We call $|V(G)|$ the order of $G$ and $| E(G )|$ the size of $G.$  For two vertices $u$ and $v$, we use the symbol $u\leftrightarrow v$ to mean that $u$ and $v$ are adjacent and use $u\nleftrightarrow v$ to mean that $u$ and $v$ are non-adjacent. Let $d_G(u,v)$ be the distance between vertices $u$ and $v$ in $G$, which represents the length of a shortest path connecting vertex $u$ and $v$ in $G$. We follow the book \cite{6} for terminology and notations.

Molecules can be modeled by graphs with vertices for atoms and edges for atomic bonds. The topological indices of a molecular graph can provide some information on the chemical properties of the corresponding molecule. It plays essential roles in the study of QSAR/QSPR in chemistry. The {\it Wiener index} of $G$, introduced in \cite{1}, is defined as
$$
 W(G)=\sum_{\{u,v\}\subseteq V(G)}d_G(u,v).
$$
Wiener \cite{1} used it to study the boiling point of paraffin. Many chemical properties of molecules are related to the Wiener index. Based on the electronic network theory, Klein and Randi\'{c} \cite{7} proposed the concept of {\it resistance distance} in 1993. The resistance distance $r_G(u,v)$ between vertices $u$ and $v$ of a connected ( molecular) graph $G$ is computed as the effective resistance between nodes $u$ and $v$ in the corresponding electrical network constructed from $G$ by replacing each edge of $G$ with an unit resistor. Without causing ambiguity, we use $r(u,v)$ for $r_G(u,v)$. This novel parameter is in fact intrinsic to the graph and has some nice interpretations and applications in chemistry (see \cite{8,9} for details). Similar to the Wiener index,  Klein and Randi\'{c} \cite{7} defined the {\it Kirchhoff index} $K\!f(G)$ of $G$ as the sum of the pairwise resistance distances between vertices, i.e.
$$
 K\!f(G)=\sum_{\{u,v\}\subseteq V(G)}r_G(u,v).
$$

Phenylenes are a class of conjugated hydrocarbons composed of $6$- and $4$-membered rings, where the $6$-membered rings (hexagons) are adjacent only to $4$-membered rings (squares), and every square is adjacent to a pair of non-adjacent hexagons. If each hexagon of a phenylene is adjacent only to two squares, we say that it is a {\it phenylene chain}. Beside benzenoid hydrocarbons, phenylenes represent another interesting class of polycyclic conjugated molecules, whose properties have been studied by a lot of researchers in the past few years. In \cite{3}, Pleter\v{s}ek  determined the edge-Wiener index and the edge-hyper-Wiener index of phenylene chain. Deng et al. \cite{4} given an efficient formula for calculating the PI index of phenylenes. Chen and Zhang \cite{10} obtained an explicit analytical expression for the Wiener index of a random phenylene chain. Li et al.\cite{2} determined the resistance distance between any two-points of the generalized phenylene. In \cite{13}, Liu and Fang obtained the explicit analytical expression for the minimal and maximal detour index of phenylene chains. Chen et al. \cite{12} determined the first five maximal (respectively,
the first five minimal) values of the Mostar index among all phenylene chains. They also characterized  the corresponding extremal chains. More interesting results for phenylene chains can be found in \cite{18}\cite{16}\cite{19}\cite{17}\cite{20}.

In this paper, we focus only on phenylene chains. Let $Q_{2n-1}$ be a ladder graph (a linear quadrilateral chain) with $2n-1$ squares. Denote by $S_i,i=1,\ldots,2n-1$ the $i$-th square of $Q_{2n-1}$. Note that a phenylene chain $G$ with $n$ hexagons and $n-1$ squares can be obtained from $Q_{2n-1}$ by add two vertices to each of the $1$-st, $3$-rd, $5$-th,\ldots, $(2n-1)$-th square. Obviously, We have three ways to add these two new vertices to $S_{2k-1},k=1,2,\ldots,n.$ That is, we can add $0$ (resp. $1$, or $2$) vertices to the top edge of $S_{2k-1}$ and the remaining vertices to the bottom edge of $S_{2k-1}.$ For convenience, we always suppose that we add two vertices to the bottom edge of $S_1$ and $S_{2n-1}$. For each of the remaining $(n-2)$-hexagons, we give a number $x_i=0$ (resp. $1$, or $2$) to the hexagon if the hexagon is obtained by adding $0$ (resp. $1$, or $2$) vertex to the top edge. In this viewpoint, we are able to represent a phenylene chain with $n$ hexagons and $n-1$ square by a $(n-2)$-vector $w=(w_1, w_2,\ldots ,w_{n-2})$ such that $w_i\in \{0,1,2\}$. In the following, we always denote a phenylene chain with $n$ hexagons and $n-1$ square by $G(w)$ such that $w$ is a $(n-2)$-tuple of $0, 1$ or $2$. A {\it kink} in a phenylene chain is a hexagon whose $w_i\neq 1.$ A phenylene chain $G(w)$ with $w_i\neq 1$ is called a {\it ``all-kink" chain}, where $1\leq i\leq n-2$. There are two special phenylene chains. We call $G(\underbrace{1,1,...,1}_{n-2})$ a {\it linear phenylene chain}, and denote it by $L_n.$ The phenylene chain $G(\underbrace{0,0,...,0}_{n-2})$ is called a {\it helicene phenylene chain}, which is denoted by $H_n$. $L_5$ and $H_5$ are illustrated in Fig. 1.
\begin{figure}[!ht]
  \centering
 \includegraphics[width=100mm]{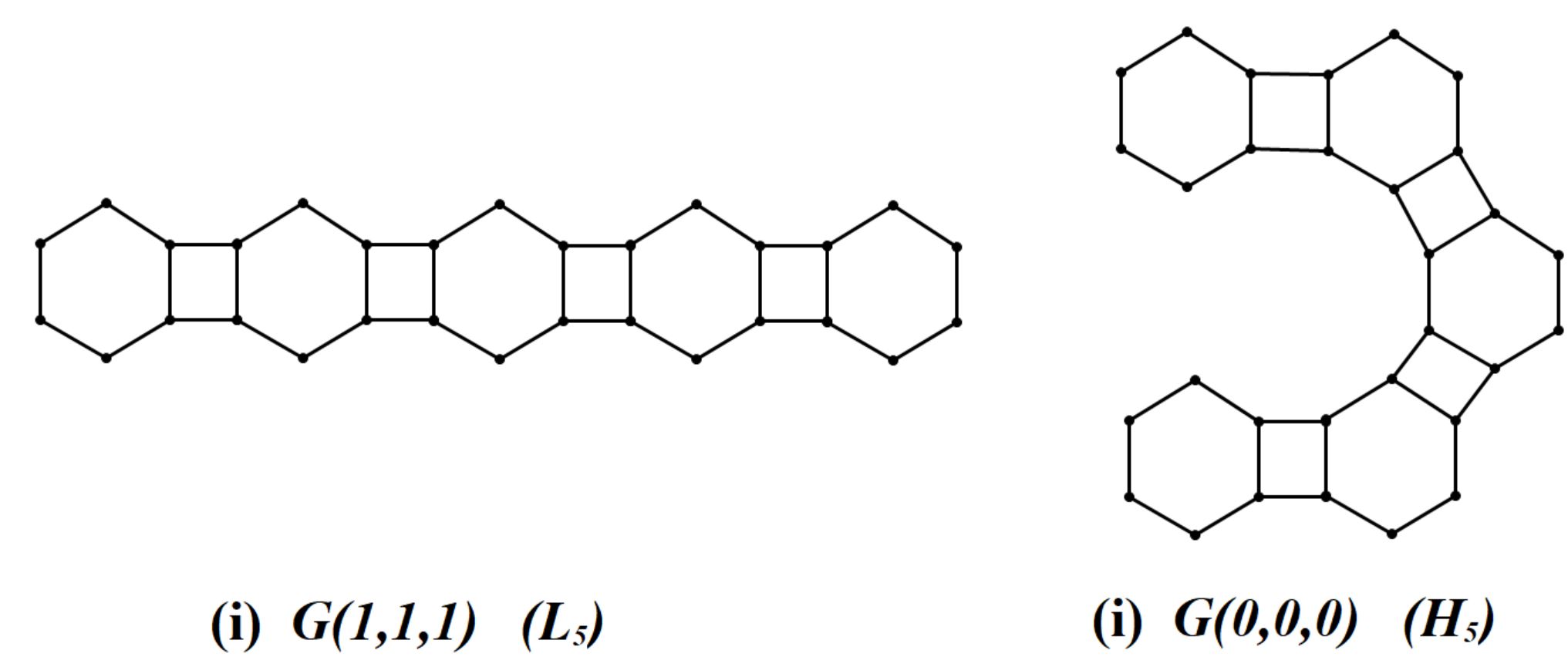}
 \caption{Some phenylene chains with $5$ hexagons and $4$ squares.}
\end{figure}

In \cite{11}, Yang and Wang proved that among all phenylene chains with $n$ hexagons and $n-1$ squares, the straight chain $L_n$ is the unique chain with maximum Kirchhoff index. In addition, they showed that among all hexagonal chains, the minimum Kirchhoff index is attained only when the phenylene chain is an ``all-kink" chain.

{\bf Theorem 1.} (\cite{11}) {\it Among all phenylene chains, the minimum Kirchhoff index is attained only when the phenylene chain is an ``all-kink" chain.}

About the minimum Kirchhoff index among all phenylene chains, Yang and Wang posed the following conjecture.

{\bf Conjecture 2.} (\cite{11}) {\it Among all phenylene chains with $n$ hexagons and $n-1$ squares, the helicene chain $H_n$ has the minimum Kirchhoff index.}

In this paper, we confirm this conjecture.

\section{Proof of the Main Result}

To prove Conjecture 2, we will need the following lemmas and definitions. There are many techniques which are used to calculate resistance distance, including the well-known series and parallel rules and the $\Delta$-$Y$ transformation, which are listed as follows.

{\bf Definition 1.} (Series Transformation) Let $x, y$, and $z$ be nodes in a graph where $y$ is adjacent to only $x$ and $z$. Moreover, let $R_1$ equal the resistance between $x$ and $y$ and $R_2$ equal the resistance between node $y$ and $z$. A series transformation transforms this graph by deleting $y$ and setting the resistance between $x$ and $z$ equal to $R_1 + R_2$.

{\bf Definition 2.} (Parallel Transformation) Let $x$ and $y$ be nodes in a multi-edged graph where $e_1$ and $e_2$ are two edges between $x$ and $y$ with resistances $R_1$ and $R_2$, respectively. A parallel transformation transforms the graph by deleting edges $e_1$ and $e_2$ and adding a new edge between $x$ and $y$ with edge resistance $r=(\frac{1}{R_1}+\frac{1}{R_2})^{-1}$

A $\Delta$-$Y$ transformation is a mathematical technique to convert resistors in a triangle formation to an equivalent system of three resistors in a $``Y"$ format as illustrated in Fig. 2. We formalize this transformation below.
\begin{figure}[!ht]
\centering
  \includegraphics[width=80mm]{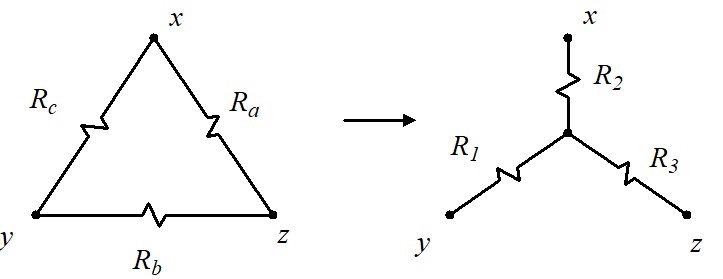}
  \caption{$\Delta$ and $Y$ circuits with vertices labeled as in Definition 3.}
\end{figure}

{\bf Definition 3.} ($\Delta$-$Y$ Transformation) Let $x,y,z$ be nodes and $R_a,R_b$ and $R_c$ be given resistances as shown in Fig. 2. The transformed circuit in the $``Y"$ format as shown in Fig. 2 has the following resistances:
\begin{eqnarray*}
  R_1=\frac{R_bR_c}{R_a+R_b+R_c}, \ \  R_2=\frac{R_aR_c}{R_a+R_b+R_c}, \ \
  R_3=\frac{R_aR_b}{R_a+R_b+R_c}.
\end{eqnarray*}

{\bf Lemma 3.} {\it \cite{5} Series transformations, parallel transformations, and $\Delta-Y$ transformations yield equivalent circuits.}

The concept of $S,T$-isomers was introduced by Polansky and Zander \cite{15} in 1982.  Suppose $A$ and $B$ are two vertex-disjoint graphs. Let $a$ and $l$ be two distinct vertices of $A$, and let $b$ and $k$ be two distinct vertices of $B$. Then $S$ is the graph obtained from $A$ and $B$ by adding edges $ab$ and $lk$. The graph $T$ is obtained from $A$ and $B$ by connecting $a$ with $k$ and $b$ with $l;$  see Figure 3.
\begin{figure}[!ht]
\centering
  \includegraphics[width=60mm]{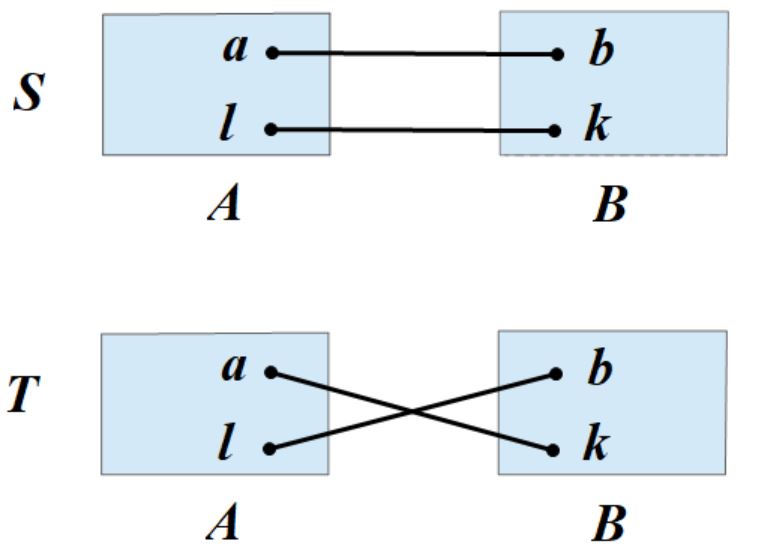}
  \caption{The structure of the graphs $S$ and $T$ and the labeling of their vertices.}
\end{figure}

Yang and Klein \cite{14} obtained the comparison theorem on the Kirchhoff index of $S,T$-isomers, which plays an essential role in the characterization of extremal phenylene chains.

{\bf Notation 1.} In the following, we will use $r_G(x)$ to denote the sum of resistance distances between $x$ and all the other vertices of $G.$ More precisely,
$$
r_G(x)=\sum_{y\in V(G)\setminus\{x\}}r_G(x,y).
$$

{\bf Lemma 4.} (\cite{14}) {\it Let $S,T,A,B,a,b,k,l$ be defined as in Fig. 3. Then}
$$
K\!f(S)-K\!f(T)=\frac{[r_A(l)-r_A(a)][r_B(b)-r_B(k)]}{r_A(a,l)+r_B(b,k)+2}.
$$

{\bf Notation 2.} Let $F$ denote a phenylene chain with $n$ hexagons and $n$ squares. For convenience, we denote the hexagons of $F$ by $C_1,C_2,\ldots,C_n$ and the squares of $F$ by $S_1,S_2,\ldots,S_{n}$ such that $C_i$ is adjacent to $S_{i}$ and $S_{i+1} \,\ (1\leq i\leq n-1)$ and $C_n$ is adjacent to $S_{n}.$  Moreover, for $1\leq i\leq n,$ we label the $i$-th 4-cycle of $F$ as $a_i,b_i,k_i,l_i.$ (See Fig. $4(i)$).
\begin{figure}[!ht]
\centering
  \includegraphics[width=90mm]{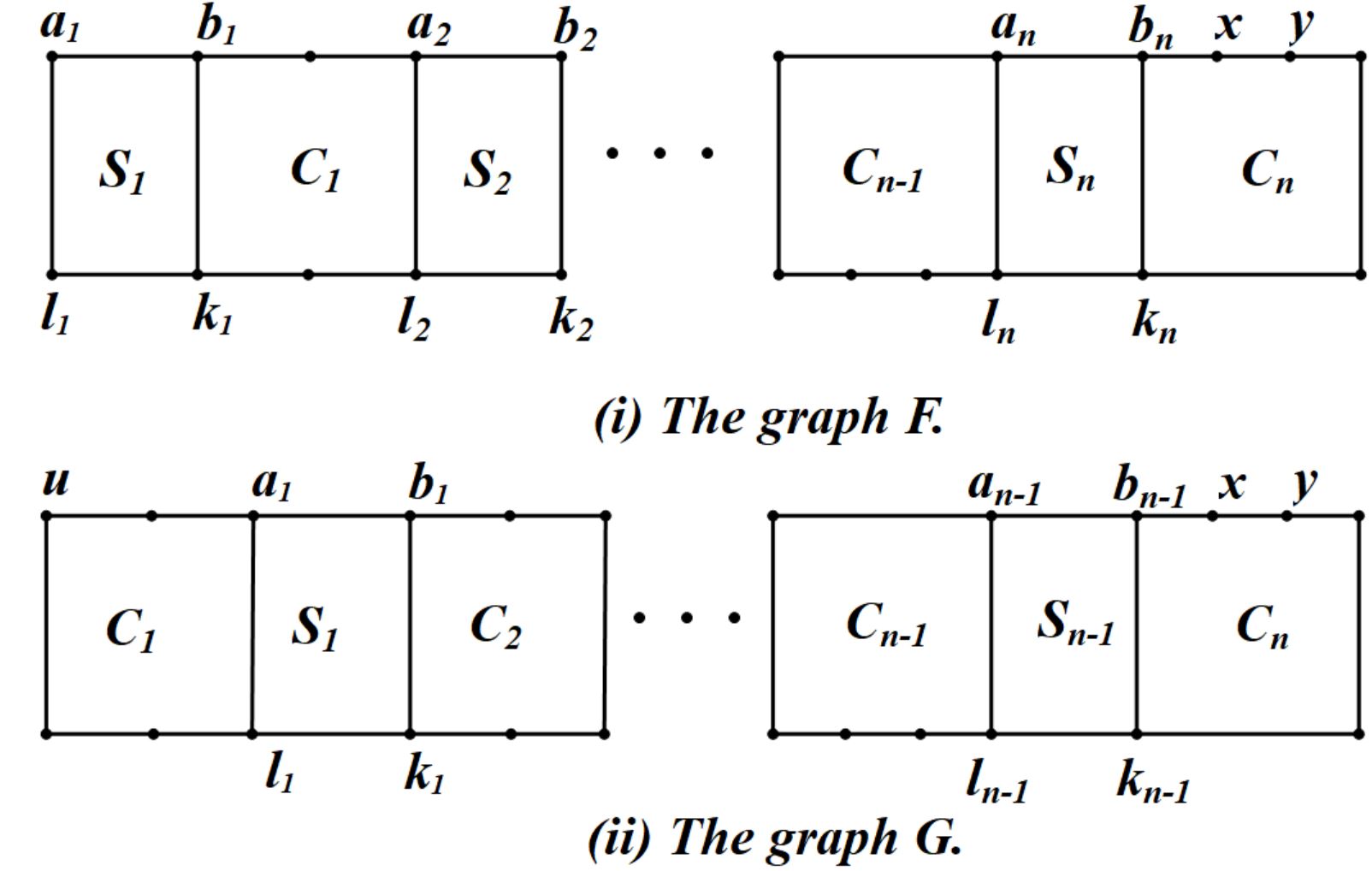}
  \caption{Figures $F$ and $G$ and their labes.}
\end{figure}
\begin{figure}[!ht]
  \includegraphics[width=140mm]{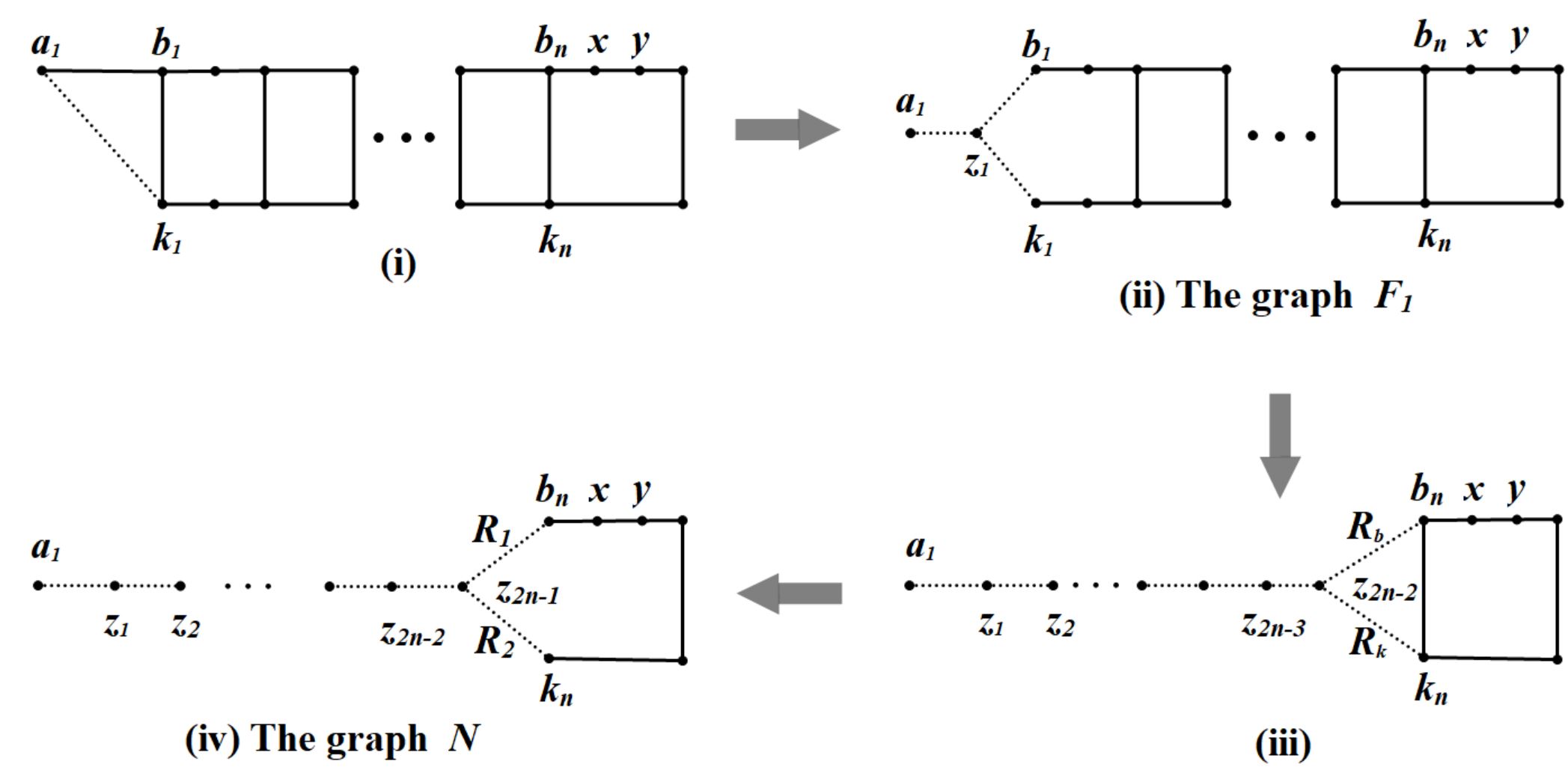}
  \caption{Illustration of circuit simplification to $F$ in the proof of Lemma 5.}
\end{figure}

{\bf Lemma 5.} {\it Let $F$ be defined as in Notation 2. Let $x$ be the vertex in $C_n$ adjacent to $b_{n}$ with degree $2,$ and let $y$ be the other neighbor of $x$. If $F$ is a weighted graph and the weight on edge $b_{n}k_{n}\in E(F)$ is $1,$ then $r_F(a_1,x)<r_F(a_1,y)$ and $r_F(l_1,x)<r_F(l_1,y).$}

{\bf Proof.} We first prove $r_F(a_1,x)<r_F(a_1,y)$. The proof of $r_F(l_1,x)<r_F(l_1,y)$ is similar to the proof of $r_F(a_1,x)<r_F(a_1,y)$, and we omit the process. In order to obtain our result, we need the following algorithm to simplify the circuit $F$.

$\bullet$ First, perform the Series transform on the leftmost cycle of $F$ to turn it into a triangle as shown in Fig. 5(i).

$\bullet$ Next, perform the $\Delta$-$Y$ transform on this new triangle. This results in a new vertex $z_1$ as shown in Fig. 5(ii). We denote this resulting graph as circuit $F_1.$

By Lemma 3, we have
$$r_F(a_1,x)=r_{F_1}(a_1,x),\,\,\,\, r_F(a_1,y)=r_{F_1}(a_1,y).$$
By repeatedly using this algorithm, we may obtain the simplified circuit of $F$ as depicted in Fig. 4(iii). Note that the edges $z_{2n-2}b_{n}$ and $z_{2n-2}k_{n}$ in Fig. $5(iii)$ are the new edges after the transformation. Denote the weighes of $z_{2n-2}b_{n}$ and $z_{2n-2}k_{n}$ by $R_b$ and $R_k.$ By making $\Delta$-$Y$ transformation to triangle $z_{2n-2}b_{n}k_n,$ we could obtain a simplified circuit $N,$ as shown in Fig. 5(iv). Suppose that the weights of edges  $z_{2n-1}b_{n}$ and $z_{2n-1}k_{n}$ in $N$ are $R_1$ and $R_2.$ Recall that the weight of edge $b_{n}k_{n}$ is $1.$ By Definition 3, we have
\begin{eqnarray*}
  R_1=\frac{R_b}{R_a+R_b+1}, \ \  R_2=\frac{R_a}{R_a+R_b+1}.
\end{eqnarray*}
Obviously, we have $0<R_1<1.$ Then using the parallel and series circuit reductions yields
\begin{eqnarray*}
  &r_F(a_1,x)=r_N(a_1,x)=r_N(a_1,z_{2n-1})+\frac{(R_1+1)(R_2+4)}{R_1+R_2+5}, \\
  &r_F(a_1,y)=r_N(a_1,y)=r_N(a_1,z_{2n-1})+\frac{(R_1+2)(R_2+3)}{R_1+R_2+5}.
\end{eqnarray*}
Hence, we have
\begin{equation*}
  r_G(a_1,x)-r_G(a_1,y)=\frac{(R_1+1)(R_2+4)}{R_1+R_2+5}-\frac{(R_1+2)(R_2+3)}{R_1+R_2+5}=\frac{R_1-R_2-2}{R_1+R_2+5}<0.
\end{equation*}
The third inequality follows from the condition $0<R_1<1.$ Therefore,we have $r_G(a_1,x)<r_G(a_1,y).$ This completes the proof. \hfill $\Box$

{\bf Notation 3.} Let $G$ denote a phenylene chain with $n$ hexagons and $n-1$ squares. For convenience, we denote the hexagons of a phenylene chain $G$ by $C_1,C_2,\ldots,C_n$ and the squares of $G$ by $S_1,S_2,\ldots,S_{n-1}$ such that $S_i$ is adjacent to $C_{i}$ and $C_{i+1} \,\ (1\leq i\leq n-1).$  Moreover, for $1\leq i\leq n-1,$ we label the $i$-th 4-cycle of $F$ as $a_i,b_i,k_i,l_i.$ (See Fig. $4(ii)$).

By a proof similar to Lemma 4, we have the following result.

{\bf Lemma 6.} {\it Let $G$ be defined as in Notation 3. Let $x$ be the vertex in $C_n$ adjacent to $b_{n-1}$ with degree $2,$ and let $y$ be the other neighbor of $x$. If $G$ is a weighted graph and the weight on edge $b_{n-1}k_{n-1}\in E(G)$ is $1,$ for any $u\in V(C_1)\setminus\{a_1,l_1\},$ we have $r_G(u,x)<r_G(u,y).$}

Now we are ready to give a proof of  Conjecture 2.
\begin{figure}[!ht]
\centering
  \includegraphics[width=100mm]{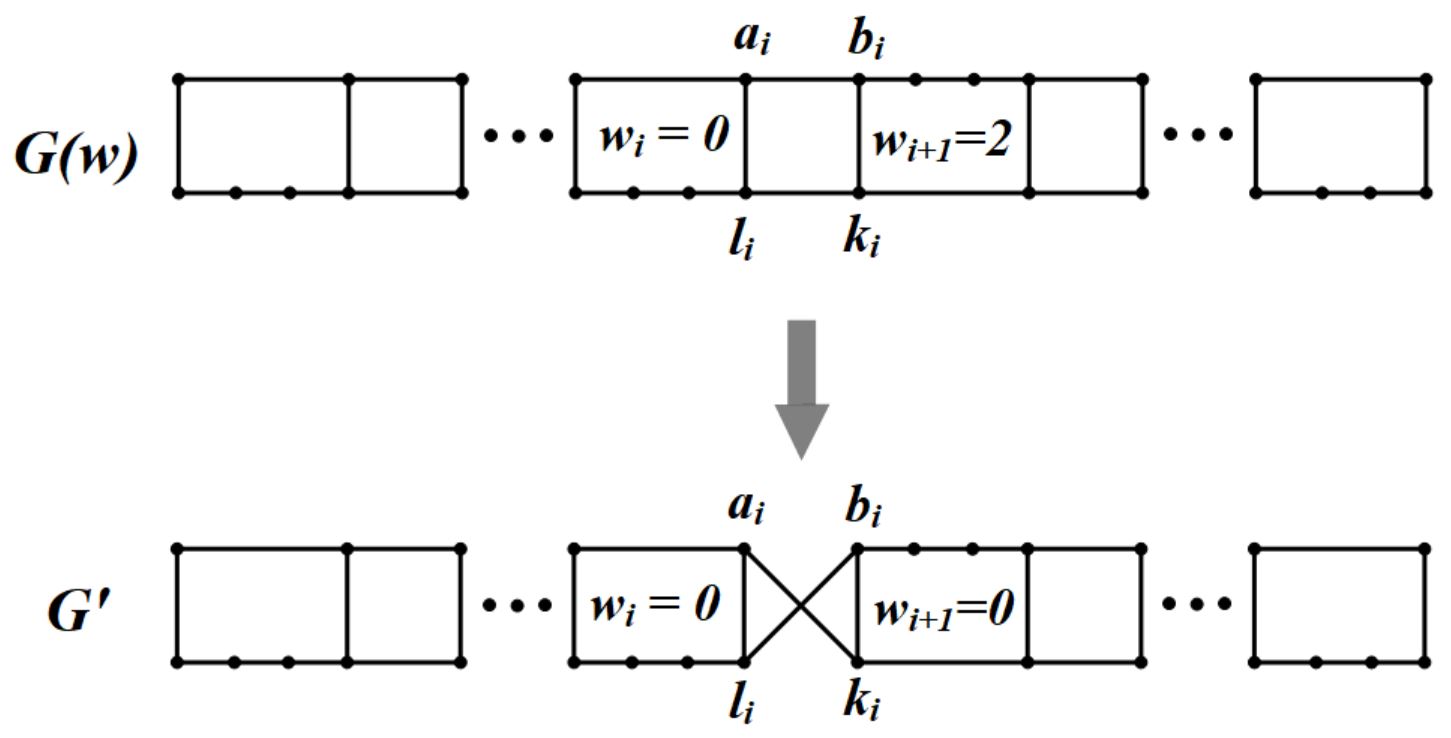}
  \caption{ $G$ and $G'$ in the proof of Conjecture 2.}
\end{figure}

{\bf Proof of Conjecture 2.} Suppose $G(w)$ is a phenylene chain and $G(w)$ has the minimum Kirchhoff index among all phenylene chains with $n$ hexagons and $n-1$ squares. Let $G(w)$ has the same labels as $G$ that be defined in Notation 3.  By Theorem 1, we have that $G(w)$ is an ``all-kink" chain. Recall that a phenylene chain $G(w)$ is an ``all-kink" chain if and only if $w$ does not contain $1.$ Without loss of generality, let $w_1=0.$ If $G(w)\neq H_n,$ there exists integer $i$ such that $w_i=0$ and $w_{i+1}=2,1\leq i\leq n-3.$ Denote the vertices of the square between $C_{i}$ annd $C_{i+1}$ by $a_i,b_i,k_i,l_i$ as show in Fig. 6. Let $G'$ be a graph obtained from $G(w)$ by deleting edge $\{a_ib_i,l_ik_i\}$ and add two new edges $\{a_ik_i,b_il_i\}.$ Next, we will prove $K\!f(G')<K\!f(G(w)).$

By the construction of $G'$, we deduce that $G(w)$ and $G'$ are pairs of $S,T$-isomers. Suppose that the two components of $G'-\{a_ik_i,b_il_i\}$ are $A$ and $B$ such that $a_i\in V(A)$ and $b_i\in V(B)$. Then by Lemma 4, we have
\begin{equation*}
K\!f(G(w))-K\!f(G')=\frac{[r_A(l_i)-r_A(a_i)][r_B(b_i)-r_B(k_i)]}{r_A(l_i,a_i)+r_B(k_i,b_i)+2}.
\end{equation*}
Now, we compare $r_A(l_i)$ and $r_A(a_i)$. Since
$$
V(A)=\left(\cup^{i-1}_{j=1} V(C_j)\setminus V(S_{j}) \right)\bigcup\left(\cup^{i-1}_{j=1} V(S_j)\setminus V(C_{j+1}) \right)\bigcup  V(C_i),
$$
we distinguish three cases. If $y\in V(C_j)\setminus V(S_{j})$, by using series and parallel connection rules, we can simplify $A$ to a weighted phenylene chain consisting of hexagons $C_j, C_{j+1},\ldots, C_{i}$ and squares $S_j, S_{j+1},\ldots, S_{i-1}.$ Note that the weight on edge $b_{i-1}k_{i-1}$ is $1.$ By Lemma 6, we have $r_A(y,a_i)<r_A(y,l_i).$ Thus,
\begin{equation}
  \sum_{y\in \cup^{i-1}_{j=1} V(C_j)\setminus V(S_{j})}r_A(y,a_i)< \sum_{y\in \cup^{i-1}_{j=1} V(C_j)\setminus V(S_{j})}r_A(y,l_i)
\end{equation}
If $y\in V(S_j)\setminus V(C_{j+1})$, by using series and parallel connection rules, we can simplify $A$ to a weighted phenylene chain consisting of squares $S_j, S_{j+1},\ldots, S_{i-1}$ and hexagons $C_{j+1}, C_{j+1},\ldots, C_{i}.$ Note that the weight on edge $b_{i-1}k_{i-1}$ is $1.$ By Lemma 5, we have $r_A(y,a_i)<r_A(y,l_i).$ Thus,
\begin{equation}
  \sum_{y\in \cup^{i-1}_{j=1} V(S_j)\setminus V(C_{j+1})}r_A(y,a_i)< \sum_{y\in \cup^{i-1}_{j=1} V(S_j)\setminus V(C_{j+1})}r_A(y,l_i)
\end{equation}
If $y\in V(C_j),$ by using series and parallel connection rules, we can simplify $A$ to a weighted hexagon $C_{i}$ with the weight $r$ on edge $b_{i-1}k_{i-1}$ and the weight $1$ on other edges.  Then
\begin{eqnarray*}
  &\sum_{y\in V(C_i)}r_A(y,a_i)=\frac{r+4}{r+5}+\frac{4(r+1)}{r+5}+\frac{3(r+2)}{r+5}+\frac{2(r+3)}{r+5}+\frac{r+4}{r+5},\\
  &\sum_{y\in V(C_i)}r_A(y,l_i)=\frac{r+4}{r+5}+\frac{2(r+3}{r+5}+\frac{3(r+2)}{r+5}+\frac{2(r+3)}{r+5}+\frac{r+4}{r+5}.
\end{eqnarray*}
Noting that the initially weight of edge $b_{i-1}k_{i-1}$ is $1,$ we have $r<1.$ Since
\begin{equation*}
  \sum_{y\in V(C_i)}r_A(y,a_i)-\sum_{y\in V(C_i)}r_A(y,l_i)=\frac{2r-2}{r+5}<0,
\end{equation*}
we have
\begin{equation}
  \sum_{y\in V(C_i)}r_A(y,a_i)<\sum_{y\in V(C_i)}r_A(y,l_i).
\end{equation}
Using Inequalities (1)-(3), we have $r_A(a_i)<r_A(l_i).$ By a similar discussion, we have $r_B(k_i)<r_B(b_i).$ Thus
\begin{equation*}
K\!f(G(w))-K\!f(G')=\frac{[r_A(l_i)-r_A(a_i)][r_B(b_i)-r_B(k_i)]}{r_A(l_i,a_i)+r_B(k_i,b_i)+2}>0.
\end{equation*}
Since $G(w)$ has the minimum Kirchhoff index, we get a contradiction.

This completes the proof. \hfill $\Box$

{\bf Acknowledgement.} The author is grateful to Prof. Xingzhi Zhan for his constant support and guidance. This research was supported by the NSFC grants 11671148 and Science and Technology Commission of Shanghai Municipality (STCSM) grant 18dz2271000.


\begin{thebibliography}{99}
\bibitem{12}H.L. Chen, H.C. Liu, Q.Q  Xiao and J.L. Zhang, Extremal phenylene chains with respect to the Mostar index. Discrete Math. Algorithms Appl. 13 (2021), no. 6, Paper No. 2150075, 27 pp.
\bibitem{10}A. Chen and F. Zhang, Wiener index and perfect matchings in random phenylene chains, MATCH Commun. Math. Comput. Chem. 61(1) (2009) 623-630.
\bibitem{18}H. Chen and Q. Guo, Tutte polynomials of alternating polyclic chains, J. Math. Chem. 57 (2019) 2248-2260.
\bibitem{16}H. Deng, J. Yang and F. Xia, A general modeling of some vertex-degree based topological indices in benzenoid systems and phenylenes, Comput. Math. Appl. 61 (2011) 3017-3023.
\bibitem{4}H. Deng, S. Chen and J. Zhang, The PI index of phenylenes, J. Math. Chem. 41 (2007) 63-69.
\bibitem{8}D.J. Klein, Resistance-distance sum rules, Croat. Chem. Acta. 75 (2002) 633-649.
\bibitem{9}D.J. Klein and O. Ivanciuc, Graph cyclicity, excess conductance, and resistance deficit, J. Math. Chem. 30 (2001) 271-287.
\bibitem{7}D.J. Klein and M. Randi\'{c}, Resistance distance, J. Math. Chem. 12 (1993) 81-95.
\bibitem{2}Q. Li, S. Li and L. Zhang, Two-point resistances in the generalized phenylenes, J. Math. Chem. 58 (2020) 1846-1873.
\bibitem{13}H.C. Liu and X.N. Fang, Extremal phenylene chains with respect to detour indices, J. Appl. Math. Comput. 67 (2021), no. 1-2, 301-316.
\bibitem{15}O.E. Polansky and M. Zander, Topological effect on MO energies, J. Mol. Struct., vol. 84, pp. 361-385, 1982.
\bibitem{3}P. Pleter\v{s}ek, The edge-Wiener index and the edge-hyper-Wiener index of phenylenes, Discrete Appl. Math. 255 (2019) 326-333.
\bibitem{19}Y. Peng and S. Li, On the Kirchhoff index and the number of spanning trees of linear phenylenes, MATCH Commun. Math. Comput. Chem. 77 (2017) 756-780.
\bibitem{5}W. Stevenson, Elements of Power System Analysis, third ed., McGraw Hill, New York, 1975.
\bibitem{6}D.B. West, Introduction to Graph Theory, Prentice Hall, Inc., 1996.
\bibitem{17}W. Wei and S. Li, Extremal phenylene chains with respect to the coefficients sum of the permanental polynomial, the spectral radius, the Hosoya index and the Merrifield-Simmons index, Discrete Appl. Math. 271 (2019) 205-217.
\bibitem{1}H. Wiener, Structural determination of paraffin boiling points, J. Am. Chem. Soc. 69 (1947) 17-20.
\bibitem{14}Y.J Yang and D.J. Klein, Comparison theorems on resistance distances and Kirchhoff indices of $S, T$-isomers, Discrete Appl. Math., vol. 175,pp. 87-93, 2014.
\bibitem{11}Y.J. Yang and D.Y. Wang, Extremal phenylene chains with respect to the Kirchhoff index and degree-based topological indices, IAENG Int. J. Appl. Math. 49 (2019), no. 3, 274-280.
\bibitem{20}Z. Zhu and J.B. Liu, The normalized Laplacian, degree-Kirchhoff index and the spanning tree numbers of generalized phenylenes, Discrete Appl. Math. 254 (2019)
256-267.
\end{thebibliography}
\end{document}